\documentclass[final,10pt,a4paper]{article}

\usepackage{graphicx}

\usepackage[T1]{fontenc}

\usepackage{amsmath}
\usepackage{amsfonts}
\usepackage{amssymb}
\usepackage{amsthm}
\usepackage{enumerate}
\usepackage{url}

\usepackage[margin=1.5in]{geometry}

\usepackage{hyperref}\hypersetup{
    colorlinks=true,
    linkcolor=blue,
    filecolor=magenta,      
    urlcolor=cyan,
    bookmarks=true,
    pdfpagemode=FullScreen,
}

\usepackage{todonotes}

\newtheorem{theorem}{Theorem}[section]
\newtheorem{lemma}[theorem]{Lemma}

\newtheorem{df}[theorem]{Definition}
\newtheorem{cor}[theorem]{Corollary}
\newtheorem{conj}[theorem]{Conjecture}
\newtheorem{assum}[theorem]{Assumption}
\newtheorem*{rem}{Remark}

\usepackage{authblk}

\title{$D(4)$-triples with two largest elements in common}
\author[a]{Marija {Bliznac Trebje\v{s}anin}}
\affil[a]{University of Split, Faculty of Science\\ Ru\dj{}era Bo\v{s}kovi\'{c}a 33, 21000 Split, Croatia\\
Email: marbli@pmfst.hr}
\date{}

\begin{document}
\maketitle

\begin{abstract} In this paper we consider two new conjectures concerning $D(4)$-quadruples and prove some special cases which support their validity. The main result is a proof that $\{a,b,c\}$ and $\{a+1,b,c\}$ cannot both be $D(4)$-triples.  
\end{abstract}

\noindent 2010 {  Mathematics Subject Classification:} 11D09, 11J68, 11J86
\\ \noindent Keywords: Diophantine $m$-tuples, Pellian equations, Hypergeometric method, Linear forms in logarithms.

\section{Introduction}
\begin{df} Let $n\neq0$ be an integer. We call a set of $m$ distinct positive integers a $D(n)$-$m$-tuple, or $m$-tuple with the property $D(4)$, if the product of any two of its distinct elements increased by $n$ is a perfect square.
\end{df}

In the classical case (when $n=1$), first studied by Diophantus, Dujella proved in \cite{duje_kon} that a $D(1)$-sextuple does not exist and that there are at most finitely many quintuples. The nonexistence of $D(1)$-quintuples was finally proven in \cite{petorke} by  He, Togb\'{e} and Ziegler.

Variants of the problem when $n=4$ or $n=-1$ are also widely studied. In the case $n=4$, similar conjectures and observations can be made as in the  case $n=1$. In the light of this observation, Filipin and the author have proven in \cite{nas2} that a $D(4)$-quintuple also does not exist. The stronger conjecture asserting the uniqueness of an extension of a triple to a quadruple with a larger element is still an open question in both cases $n=1$ and $n=4$. Moreover, in the case $n=-1$, a conjecture about the nonexistence of a quadruple is studied. 

Let $\{a,b,c\}$, $a<b<c$, be a $D(4)$-triple. We define
$$d_{\pm}=d_{\pm}(a,b,c)=a+b+c+\frac{1}{2}\left(abc\pm \sqrt{(ab+4)(ac+4)(bc+4)}\right).$$
It is straightforward to check that $\{a,b,c,d_{+}\}$  is a $D(4)$-quadruple, which we will call a regular quadruple.  If $d_{-}\neq 0$ then $\{a,b,c,d_{-}\}$ is also a regular $D(4)$-quadruple with $d_{-}<c$ and $c=d_+(a,b,d_{-})$. In other words, it is conjectured that an irregular $D(4)$-quadruple doesn't exist.

Results that support this conjecture in some special cases can be found for example in \cite{aft}, \cite{bf}, \cite{dujram}, \cite{fil_par}, \cite{fht}. 
In \cite{mbt}, the author has proved that a $D(4)$-quadruple $\{a,b,c,d\}$ with $a<b<c<d$ and $c\geq  39247b^4$ must be a regular $D(4)$-quadruple.

\medskip
 Let us describe a problem of extension of a $D(4)$-triple $\{a,b,c\}$ to a quadruple with an element $d$. Then the element $d$ satisfy equalities
$$ad+4=x^2,\quad bd+4=y^2, \quad cd+4=z^2,$$
where $x,y,z$ are some positive integers.
A system of generalized Pellian equations
\begin{align}
cx^2-az^2&=4(c-a), \label{pellova_cetv_ac}\\
cy^2-bz^2&=4(c-b), \label{pellova_cetv_bc}
\end{align}
is obtained by eliminating $d$ from previous equations. It is not hard to describe the sets of solutions of equations (\ref{pellova_cetv_ac}) and (\ref{pellova_cetv_bc}), see for example \cite{fil_xy4}. The main approach in solving this problem is finding an upper bound for the number $z$ using the hypergeometric method and Baker's method.

 So far, research has shown that variants of the problem with $n=1$ and $n=4$ are closely related in results and methods used to prove them but differ in the details of the proof, which will also be the case here. This paper will closely follow ideas and methods from \cite{cdf} to prove analogous results in the case of $D(4)$-$m$-tuples.  

As in \cite{cdf}, we expect that the following conjectures hold and prove two theorems that support their validity.

\begin{conj}\label{conj_1}
Suppose that $\{a_1,b,c\}$ and $\{a_2,b,c\}$ are $D(4)$-triples with $a_1<a_2<b<c$. Then, $\{a_1,a_2,b,c\}$ is a $D(4)$-quadruple.
\end{conj}

Since it has been proved in \cite{nas2} that a $D(4)$-quintuple cannot exist, the next conjecture also follows from the previous one. 

\begin{conj}\label{conj_2}
Suppose that $\{a_1,b,c,d\}$ is a Diophantine quadruple with $a_1<b<c<d$. Then, $\{a_2,b,c,d\}$ is not a Diophantine quadruple for any integer $a_2$ with $a_1\neq a_2<b$.
\end{conj} 

Conjecture \ref{conj_1}\, asserts if $\{a_1,a_2\}$ is not a $D(4)$-pair then $\{a_1,b,c\}$ and $\{a_2,b,c\}$ cannot both be $D(4)$-triples. In the next theorem, we will observe some pairs of a form $\{a,a+1\}$ and prove that they have that desired property, which supports the claim of Conjecture \ref{conj_1}.\par
The only $\{a,a+1\}$ $D(4)$-pair is $\{3,4\}$ and it can be extended to infinitely many different quadruples $\{3,4,c,d_+\}$ (they are explicitly described in \cite{bf}). For example, one of them is $\{3,4,15,224\}$, so $\{3,15,224\}$ and $\{4,15,224\}$ are $D(4)$-triples. We will prove that the same cannot hold for any other positive integer $a$. 


\begin{theorem} \label{tm_ap1}
Suppose that $\{a,b,c\}$ is a $D(4)$-triple, $a\neq 3$. Then, $\{a+1,b,c\}$ is not a $D(4)$-triple. 
\end{theorem}

Proof of this theorem will be separated in two cases in Section \ref{dokaz tm}, the first case will be proven by using the hypergeometric method and the second case by using linear forms in logarithms.
  
We further support the validity of conjectures by proving the next results.
\begin{theorem}\label{tm_conjec}
   If $c<0.25b^3$, then Conjecture \ref{conj_1} holds.
\end{theorem}
As a consequence of Theorem 1.6 from \cite{mbt} and Theorem 1 from \cite{nas2}, one sees that Conjecture \ref{conj_2} holds when $c\geq 39247b^4$. 
\begin{cor}
   If either $c<0.25b^3$ or $c\geq 39247b^4$, then Conjecture \ref{conj_2} holds. 
\end{cor}

\section{Pellian equations and preliminary results}

Let $\{a,b,c\}$ be a $D(4)$-triple, $a\neq 3$. Suppose $\{a+1,b,c\}$ is also a $D(4)$-triple. Without loss of generality suppose $b<c$. There exits positive integers $s,t$ such that 
\begin{align*}
    ab+4&=s^2\\
    (a+1)b+4&=t^2.
\end{align*}
We get a Pellian equation
\begin{equation}\label{pell_3}
at^2-(a+1)s^2=-4,
\end{equation}
with solutions for unknown $s$ given by a recurrent sequence
\begin{equation}\label{sovi}
s_0=2,\quad s_1=8a+2,\quad s_{\nu+2}=2(2a+1)s_{\nu+1}-s_{\nu},\quad \nu\in\mathbb{N}_0.
\end{equation}

Define
$b_{\nu}=(s_{\nu}^2-4)/a$. Then we can express element $b$ in the terms of $a$ and use these values in our proof. Explicitly,
\begin{align*}
b_1=&64a+32,\\ 
b_2=&1024a^3+1536a^2+704a+96,\\
b_3=&16384a^5+40960a^4+37888a^3+15872a^2+2944a+192,\\
b_4=&262144a^7+917504a^6+1294336a^5+942080a^4+375808a^3+80384a^2\\
&\quad+80384a+320.
\end{align*}

Also, there exist positive integers $x,y,z$ such that 
\begin{align*}
    ac+4&=x^2,\\
    (a+1)c+4&=y^2,\\
    bc+4&=z^2.
\end{align*}
These equations give a system of Pellian equations
\begin{align}
    az^2-bx^2&=4(a-b), \label{prva_pelova_s_a}\\
    (a+1)z^2-by^2&=4(a+1-b),\label{druga_pelova_s_ap1}
\end{align}
whose solutions $(z,x)$ and $(z,y)$ we will further observe. As in \cite[Lemma 2]{dujram} we can describe solutions of this system. 
\begin{lemma}\label{granice_fundamentalnih}
Let $(z,x)$ and $(z,y)$ be positive solutions of (\ref{prva_pelova_s_a}) and (\ref{druga_pelova_s_ap1}). Then there exist solutions $(z_0,x_0)$ of (\ref{prva_pelova_s_a}) and $(z_1,y_1)$ of (\ref{druga_pelova_s_ap1}) in the ranges
\begin{align*}
1&\leq x_0<\sqrt{\frac{a(b-a)}{s-2}},\\
1&\leq |z_0|<\sqrt{\frac{(s-2)(b-a)}{a}},\\
1&\leq y_1<\sqrt{\frac{(a+1)(b-a-1)}{t-2}},\\
1&\leq |z_1|<\sqrt{\frac{(t-2)(b-a-1)}{a+1}},
\end{align*}
such that 
\begin{align*}
z\sqrt{a}+x\sqrt{b}&=(z_0\sqrt{a}+x_0\sqrt{b})\left(\frac{s+\sqrt{ab}}{2}\right)^m,\\
z\sqrt{a+1}+y\sqrt{b}&=(z_1\sqrt{a+1}+y_1\sqrt{b})\left(\frac{t+\sqrt{(a+1)b}}{2}\right)^n.
\end{align*}
\end{lemma}
As before, solutions can be expressed as elements of recurrent sequences. More precisely, $z$ must be an element of sequences
\begin{align}
&v_0=z_0,\ v_1=\frac{1}{2}\left(sz_0+bx_0\right),\ v_{m+2}=sv_{m+1}-v_{m},\label{vm}\\
&w_0=z_1,\ w_1=\frac{1}{2}\left(tz_1+by_1 \right),\ w_{n+2}=tw_{n+1}-w_n,\label{wn}
\end{align}
where $m,n\geq 0$ are positive integers.
Following this notation, it must hold $z=v_m=w_n$ for some $m$ and $n$.

It is easy to see that
\begin{equation}z_0^2\equiv z_1^2\equiv 4\ (\bmod \ b).\label{kong_4}\end{equation} 

For simplicity of the proof, we will assume that $b$ and $c$ are "minimal" among all $b$'s and $c$'s satisfying the conditions of Theorem \ref{tm_ap1}. 
\begin{assum}\label{pretpostavka}
At least one of $\{a,b',b\}$ and $\{a+1,b',b\}$ is not a $D(4)$-triple for any $b'$ with $0<b'<b$.
\end{assum}

Since we are searching for intersections of sequences (\ref{vm}) and (\ref{wn}), we can describe the initial terms of sequences more precisely. We omit the proof since it is proven similarly as \cite[Lemma 2.3]{cdf} by following cases from \cite[Lemma 4]{fil_xy4}.
\begin{lemma}\label{lema_fund_parno}
If the equation $v_m=w_n$ has a solution, then both $m$ and $n$ are even and $z_0=z_1=2\varepsilon$, where $\varepsilon\in\{\pm 1\}$.
\end{lemma}

\begin{rem}
Assumption \ref{pretpostavka} is crucial for the proof of Lemma \ref{lema_fund_parno}.
If $a=3$, then (\ref{pell_3})  would also have a fundamental solution $(t_0,s_0)=(0,1)$. In this case, for every $b=b_{\nu}=(s_{\nu}^2-4)/a,$ where $s_{\nu}$ is defined as in (\ref{sovi}), there would exist $b'<b$ which arises from the sequence with this second fundamental solution. 
\end{rem}

From $z_0=z_1=2\varepsilon$ we have $x_0=y_1=2$ so sequences (\ref{vm}) and (\ref{wn}) can be written in the form
\begin{align}
    &v_0=2\varepsilon,\ v_1=\varepsilon s+b,\ v_{m+2}=sv_{m+1}-v_{m},\label{vm_epsilon}\\
&w_0=2\varepsilon,\ w_1=\varepsilon t+b,\ w_{n+2}=tw_{n+1}-w_n.\label{wn_epsilon}
\end{align}

Following lemma is easily proved by induction.

\begin{lemma} \label{lema_kongruencije}
\begin{align*}
    v_{2m}&\equiv 2\varepsilon+b(a\varepsilon m^2+sm)\ (\bmod \ b^2),\\
    v_{2m+1}&\equiv \varepsilon s+b\left(\frac{1}{2}as\varepsilon m(m+1)+(2m+1)\right)\ (\bmod \ b^2),\\
    w_{2n}&\equiv 2\varepsilon+b((a+1)\varepsilon n^2+tn)\ (\bmod \ b^2),\\
    w_{2n+1}&\equiv \varepsilon t+b\left(\frac{1}{2}(a+1)t\varepsilon n(n+1)+(2n+1)\right)\ (\bmod \ b^2).
\end{align*}
\end{lemma}

By using previous results we can give some upper and lower bounds on indices $m$ and $n$. 
\begin{lemma}\label{lema_m_n_odnos}
If $v_m=w_n$ has a solution, then $n\leq m\leq \frac{3}{2}n+1$.
\end{lemma}
\begin{proof}
 We follow the proof of \cite[Lemma 5]{fil_xy4}. By using that $b\geq b_1\geq 96a\geq  96$
 we get the estimates for $v_m$ and $w_n$ which yields an inequality
 $$(s-1)^{m-1}<2.396t^n<t^{n+0.332}.$$
 Also, $(\sqrt{ab+4}-1)^3=(s-1)^3>t^2=(a+1)b+4$ so we conclude 
 $$m\leq \frac{3}{2}n+1.$$
 Similarly, by observing
 $$(t-1)^{n-1}<s^{m+0.451}<(t-1)^{m+0.451},$$
 we get $n-1\leq m$. Now the desired inequality follows from Lemma \ref{lema_fund_parno}.
\end{proof}

Observe that $m=0$ means that $z=v_0=2\varepsilon $ and  $bc+4=4$, i.e., $bc=0$ which cannot hold. So, Lemma \ref{lema_fund_parno} implies $m\geq 2$.
\begin{lemma}
If $v_m=w_n$ has a solution with $m\geq 2$ then $m>n$.
\end{lemma}
\begin{proof}
We will prove by induction that $v_n\leq w_n$ for $n\geq 2$, so if $v_m=w_n$ then $m>n$ for $m\geq 2$.\par
Let $n=2$. We first observe the case  $\varepsilon=1$, where $x_0=y_1=2$ so $v_2=b(a+s)+2<b(a+1+t)+2=w_2$. On the other hand, if $\varepsilon=-1$, first we observe that $s+1<t$, since $s+1\geq t$ is in a contradiction with $b\geq b_1=64a+32$. Using that inequality, we have $v_2=b(s-a)-2<b(t-1-a)-2=w_2$.\\
Let us assume that $v_{n-1}<w_{n-1}$ for some $n\geq 3$. Since $s\leq t-2$ and sequence $w_n$ is increasing, we have
$$v_n=sv_{n-1}-v_{n-2}<sv_{n-1}<sw_{n-1}\leq tw_{n-1}-2w_{n-1}<tw_{n-1}-w_{n-2}=w_n,$$
which proves our statement. 
\end{proof}


\begin{lemma}\label{donja_za_m}
If $v_m=w_n$ has a solution with $m>0$ then 
$$m>0.4672(a+1)^{-1/2}b^{1/2}.$$
\end{lemma}
\begin{proof}
Since $v_m=w_n$, from Lemma \ref{lema_kongruencije} we get that congruence
\begin{equation}\label{kongruencija}
\varepsilon(am^2-(a+1)n^2)\equiv tn-sm\ (\bmod \ b)
\end{equation}
holds. 

Suppose to the contrary that $m\leq 0.4672(a+1)^{-1/2}b^{1/2}$. Observe that
$$\max\{am^2,(a+1)n^2\}\leq (a+1)m^2<0.4672^2b<0.5b$$
and 
$$\max\{sm,tn\}\leq tm<0.4672\sqrt{1+\frac{4}{(a+1)b}}\cdot b<0.4672\sqrt{\frac{49}{48}}b<0.5b,$$
since $(a+1)b\geq (a+1)b_1\geq 192$. This implies that congruence (\ref{kongruencija}) is an equality. After multiplying with $tn+sm$ we get
\begin{equation}\label{jdba2}
    ((a+1)n^2-am^2)(b+\varepsilon(tn+sm))=4(m^2-n^2).
\end{equation}
Since $m>n$, both side of equation are positive. Also, by using Lemmas \ref{lema_fund_parno} and \ref{lema_m_n_odnos} on equation (\ref{jdba2}) we get
$$|b+\varepsilon(tn+sm)|\leq m^2-n^2$$
must hold. Then
\begin{align*}
    b&\leq tn+sm+m^2+n^2\leq t(m-2)+sm+\frac{5}{9}m^2=m\left(t\left(1-\frac{2}{m}\right)+s+\frac{5}{9}m\right)\\
    &\leq b\cdot 0.4672\left[\sqrt{1+\frac{4}{b(a+1)}}\left(1-\frac{2}{0.4672}\sqrt{\frac{a+1}{b}}\right)+\sqrt{\frac{ab+4}{b(a+1)}}+\frac{5}{9}\cdot\frac{0.4672}{a+1}\right].
\end{align*}
Since $a\geq 1$ and $b\geq b_1=64a+32\geq 96$ we get
$$b<b\cdot 0.4672\left(\sqrt{\frac{49}{48}}\cdot 1+1+\frac{5}{18}\cdot 0.4672\right)<b,$$
a contradiction. 
\end{proof}

Elements of the sequences $(v_m)$ and $(w_m)$ can be expressed explicitly 
\begin{align}
&v_m=\frac{\varepsilon \sqrt{a}+\sqrt{b}}{\sqrt{a}}\left(\frac{s+\sqrt{ab}}{2}\right)^m+\frac{\varepsilon \sqrt{a}-\sqrt{b}}{{\sqrt{a}}}\left(\frac{s-\sqrt{ab}}{2}\right)^m,\label{vm_eksplicitno}\\
&w_n=\frac{\varepsilon \sqrt{a+1}+\sqrt{b}}{\sqrt{a+1}}\left(\frac{t+\sqrt{(a+1)b}}{2}\right)^n+\frac{\varepsilon \sqrt{a+1}-\sqrt{b}}{\sqrt{a+1}}\left(\frac{t-\sqrt{(a+1)b}}{2}\right)^n.\label{wn_eksplicitno}\end{align}

If $z=v_m=w_n$ we define a linear form in three logarithms 
$$\Lambda:=m\log\alpha-n\log\beta+\log\gamma,$$
where $$\alpha=\frac{s+\sqrt{ab}}{2},\  \beta=\frac{t+\sqrt{(a+1)b}}{2} \textrm{ and } \gamma=\frac{\sqrt{a+1}(\sqrt{b}+\varepsilon\sqrt{a})}{\sqrt{a}(\sqrt{b}+\varepsilon\sqrt{a+1})}.$$
It is not hard to show that (see \cite[Lemma 10]{fil_xy4})
\begin{equation}\label{lambda}
0<\Lambda<\alpha^{1-2m}.
\end{equation}
By using the fact that $b\geq b_1\geq 64a$ and (\ref{vm_eksplicitno}) it is easy to see that the next lemma holds.
\begin{lemma}\label{lema_vm_za_linearne}
If $m\geq 1$, then $z=v_m>\left(\dfrac{s+\sqrt{ab}}{2}\right)^m.$
\end{lemma}
The next lemma can be proven by following the idea of \cite[Lemma 4.1]{cdf}.
\begin{lemma}\label{lema_m_b1}
If $v_m=w_n$ has a solution with $m\geq 2$, then
$$(m-0.0005)\log\alpha-n\log\beta<0.$$
\end{lemma}

\section{Proof of Theorem \ref{tm_ap1}}\label{dokaz tm}
The next theorem is part of the hypergeometric method first developed in \cite{rickert}. We omit the details of its proof since it only slightly differs from the proof of \cite[Theorem 3.2]{cdf} or \cite[Theorem 2]{bf}.
\begin{theorem}\label{tm_rickert}
Let $a$ be a positive integer and $N$ a multiple of $a(a+1).$ Assume that $N\geq270a(a+1)^2$. Then the numbers $\theta_1=\sqrt{1+4a/N}$ and $\theta_2=\sqrt{1+4(a+1)/N}$ satisfy
$$\max\left\{\left|\theta_1-\frac{p_1}{q}\right|,\left|\theta_2-\frac{p_2}{q}\right|\right\}>(2.96\cdot 10^{28}N(a+1))^{-1}q^{-\lambda}$$
for all integers $p_1,p_2,q$ with $q>0$, where
$$\lambda=1+\frac{\log(11(a+1)N)}{\log(0.041a^{-1}(a+1)^{-1}N^2)}<2.$$
\end{theorem}

\begin{lemma}[cf. {\cite[Lemma 14]{fil_xy4}}]\label{lema_z}
Let $N=a(a+1)b$ and let $\theta_1$, $\theta_2$ be as in Theorem \ref{tm_rickert}. Then, all positive solutions to the system of Pellian equations (\ref{prva_pelova_s_a}) and (\ref{druga_pelova_s_ap1}) satisfy
$$\max\left\{\left|\theta_1-\frac{(a+1)sx}{a(a+1)z}\right|,\left|\theta_2-\frac{aty}{a(a+1)z}\right|\right\}<\frac{2b}{a}z^{-2}.$$
\end{lemma}

\begin{proof}(\textbf{Proof of the Theorem \ref{tm_ap1} in the case $b\geq b_2$})\\

Let us assume that $b\geq b_2=1024a^3+1536a^2+704a+96$ and that Assumption \ref{pretpostavka} holds. We can apply Theorem \ref{tm_rickert} with $N=a(a+1)b$, $p_1=(a+1)sx$, $p_2=aty$ i $q=a(a+1)z$ and Lemma \ref{lema_z} to show that the next inequality holds
$$z^{2-\lambda}<5.92\cdot 10^{28}(a+1)^2b^2a^{\lambda}(a+1)^{\lambda}.$$
Inserting expression for $\lambda$ and approximating $\lambda<2$ on the right hand side of the previous inequality yields
\begin{equation}\label{njd_z}
\log z<\frac{\log(5.92\cdot 10^{28}a^2(a+1)^4b^2)\log(0.041a(a+1)b^2)}{\log(0.0037(a+1)^{-1}b)}.
\end{equation}

Combining Lemmas \ref{donja_za_m} and \ref{lema_vm_za_linearne} with inequality (\ref{njd_z})  implies
\begin{equation}\label{njd_treca}0.4672(a+1)^{-1/2}b^{1/2}<\frac{\log(5.92\cdot 10^{28}a^2(a+1)^4b^2)\log(0.041a(a+1)b^2)}{\log\left(\frac{s+\sqrt{ab}}{2}\right)\log(0.0037(a+1)^{-1}b)}.\end{equation}
Since $b\geq b_2\geq1024a^2(a+1)>1024a^3$ and the right-hand side of the inequality is decreasing in $b$, leaving us to observe an inequality
$$14.95a<\frac{\log(6.21\cdot 10^{34}a^6(a+1)^6)\log(42992a^5(a+1)^3)}{\log\left(32a^2\right)\log(3.78a^2)}.$$
If $a>5$ (i.e. $a\geq 6$), applying $a+1<1.2a$ to the previous inequality and solving for $a$ returns  $a\leq 5$, a contradiction. 
On the other hand, for $a\leq 5$ we can insert values for $a$ and $b_3$ in the inequality (\ref{njd_treca}) and see it cannot hold, which means it remains to consider only the pairs $(a,b_2)$,  $1\leq a\leq 5$.

 Since $b$ achieves its maximum for $a=5$, we can use it in the inequality from the proof of \cite[Theorem 1]{sestorka}
$$\frac{m}{\log(m+1)}<6.543\cdot 10^{15}\log^2 b$$
to get $m<4.3\cdot 10^{19}$. Also, from Lemma \ref{donja_za_m} and minimal value for $b$ we have $m\geq 0.4672\sqrt{b/5}>12$. 
Now we can use the Baker-Davenport reduction method on the remaining pairs as described in \cite{dujpet}, and each case returned the bound $m<7$, which is a contradiction. 
\end{proof}
\medskip 

It remains to consider the case $b=b_1$.  Lemma \ref{lema_m_b1} can be easily applied to get another useful relation between indices $m$ and $n$. 
\begin{lemma}\label{lema_n_m_b1}
If $v_m=w_n$ has a solution with $m\geq 2$ and $b=b_1=64a+32$, then
$$n>2(\nu-0.0005)a\log \alpha,$$
where $\nu=m-n$.
\end{lemma}

\begin{theorem}(\cite[Corollary 2]{lm})\label{tm_lm} Assume that $\alpha_1$ and $\alpha_2$ are real, positive and multiplicatively independent algebraic numbers in a field $K$ of degree $D$. Set
$$\Lambda:=b_2\log \alpha_2-b_1\log\alpha_1,$$
where $b_1$ and $b_2$ are positive integers. Let $A_1$ and $A_2$ be real numbers greater than one such that 
$$\log A_i\geq \max\{h(\alpha_i),|\log\alpha_i|/D,1/D\},\quad (i=1,2).$$
Set
$$b':=\frac{b_1}{D\log A_2}+\frac{b_2}{D\log A_1}.$$
Then,
$$\log\Lambda>-24.34D^4(\max\{\log b'+0.14,21/D,1/2\})^2\log A_1\log A_2.$$

\end{theorem}

\begin{proof}(\textbf{Proof of the Theorem \ref{tm_ap1} in the case $b= b_1$})\\
Let us denote $\nu=m-n$ and rewrite $\Lambda $ as follows
$$\Lambda=\log(\alpha^{\nu}\gamma)-n\log(\beta/\alpha).$$
Let
$$b_1=n,\ b_2=1,\ \alpha_1=\beta/\alpha,\ \alpha_2=\alpha^{\nu}\gamma,\ D=4.$$
Multiplicative independence of $\alpha_1$ and $\alpha_2$ over $\mathbb{Q}(\sqrt{ab},\sqrt{(a+1)b})$ can be verified similarly as in \cite[Lemma 19]{petorke}, so the linear form in logarithms $\Lambda$ with these parameters satisfy conditions of Theorem \ref{tm_lm}.

We have $$h(\alpha)=\frac{1}{2}\log\alpha,\quad h(\beta)=\frac{1}{2}\log \beta.$$
By observing conjugates of $\gamma$ whose absolute values are greater than one (see \cite{cdf}) and noting that the leading coefficient of the minimal polynomial of $\gamma$ is divisor of $a^2(b-a-1)^2$ we estimate 
$$h(\gamma)\leq \frac{1}{4}\log\left[a^{1/2}(a+1)^{3/2}(b-a)(\sqrt{b}+\sqrt{a})(\sqrt{b}+\sqrt{a+1})\right]<\log(2\alpha).$$
This implies
\begin{align*}
    h(\alpha_1)&=h(\beta/\alpha)\leq h(\beta)+h(\alpha)=\frac{1}{2}(\log \alpha+\log \beta),\\
    h(\alpha_2)&=h(\alpha^{\nu}\gamma)\leq \nu h(\alpha)+h(\gamma)\leq \left(\frac{\nu}{2}+1\right)\log (\alpha)+\log 2<1.16\left(\frac{\nu}{2}+1\right)\log (\alpha),
\end{align*}
where we have used that $ab\geq 96$ and $\nu\geq 2$.
Since $\gamma\leq 2$, we also have   $$\frac{\log \alpha_2}{D}<\left(\frac{\nu}{2}+1\right)\log (\alpha)<1.16\left(\frac{\nu}{2}+1\right)\log (\alpha),$$
hence we may take
\begin{align*}
\log A_1&=\frac{1}{2}(\log\alpha+\log\beta)\\
\log A_2&=1.16\left(\frac{\nu}{2}+1\right)\log\alpha.
\end{align*}
Now
$$b'=\frac{n}{2\cdot 1.16\cdot (\nu+2)\log\alpha}+\frac{1}{2(\log\alpha+\log\beta)}<\frac{m}{2(\nu+2)\log\alpha}.$$
By using $b\geq 96$ we can estimate that
$\beta<1.43\alpha$, hence $\log\alpha+\log\beta<\log(1.43\alpha^2)<2.16\log\alpha$.
From Theorem \ref{tm_lm} and (\ref{lambda}) it follows
\begin{equation}
    \label{njdn_4}
    \frac{1.16m-0.58}{2(\nu+2)\log\alpha}<1132\left(\max\left\{\log\frac{1.16m}{2(\nu+2)\log\alpha},5.25\right\}\right)^2.
\end{equation}
If $\log\frac{1.16m}{2(\nu+2)\log\alpha}>5.25$, the inequality (\ref{njdn_4}) implies that
\begin{equation}\label{gornja_m}
m<18067.6(\nu+2)\log \alpha.
\end{equation}
In the other case, when $\log\frac{1.16m}{2(\nu+2)\log\alpha}\leq 5.25$ the same inequality holds. By combining Lemma \ref{lema_n_m_b1} and (\ref{gornja_m}) we obtain
$$a<9033.8\frac{\nu+2}{\nu-0.0005}<18072.11.$$
Therefore it remains to verify only finitely many pairs ${a,b_1}$, $1\leq a\leq 18072$. Three steps of the Baker-Davenport reduction method ended with the bound $m<2$, which is a contradiction.
\end{proof}

\section{Proof of Theorem \ref{tm_conjec}}
Let $\{a_1,b,c\}$ and $\{a_2,b,c\}$ be $D(4)$-triples with $a_1<a_2<b<c$ and $c<0.25b^3$. Define
$$d_i=a_i+b+c+\frac{1}{2}(a_ibc-r_is_iu),$$
where $r_i=\sqrt{a_ib+4}$, $s_i=\sqrt{a_ic+4}$ and $u=\sqrt{bc+4}$.
It is easy to see, if $d_i>0$ then $\{a_i,d_i,b,c\}$ is also a $D(4)$-quadruple with $d_i< c$ and $c=d_+(a,b,d_i)$. There is also a possibility that $d_i=0$, namely in the case $c=a_i+b+2r$. In both cases the relation 
\begin{equation}\label{relacija}
(b+c-a_i-d_i)^2=(a_id_i+4)(bc+4)
\end{equation}
holds. We will denote $t_i=\sqrt{a_id_i+4}$. 

There exist a rational number $\lambda_i$ satisfying
$$c=abd_i+\lambda_i\max\{d_i,b\},\quad 1<\lambda_i<4,$$
implying 
\begin{equation}\label{odnos_di}
    (a_1d_1-a_2d_2)b=\lambda_2\max\{d_2,b\}-\lambda_1\max\{d_1,b\}.
\end{equation}
There are three cases to consider depending on the value of $\max\{b,d_1,d_2\}$. In each case we obtain that  $a_1d_1=a_2d_2$ must hold by following similar arguments, so we give details only for the case $d_1=\max\{b,d_1,d_2\}$. Assume on the contrary, that $a_1d_1\neq a_2d_2$.  Then
$|a_1d_1-a_2d_2|=|t_1^2-t_2^2|\geq 2t_1-1,$
which together with (\ref{odnos_di}) implies 
\begin{equation}\label{gornja_t}2t_1<4\frac{d_1}{b}+1.\end{equation} 
After squaring we get an inequality $d_1>\frac{1}{4}b(a_1b-2)$ and inserting it in the definition for $t_1$ we get 
$$t_1 \geq \frac{1}{2}a_1b.$$
Combining with (\ref{gornja_t}) implies a slightly better lower bound $d_1>\frac{1}{4}b(a_1b-1)$ which finally gives
$$c=c_+(a_1,b,d_1)>b+d_1+a_1d_1b>b+\frac{1}{4}b(a_1b-1)+\frac{1}{4}a_1b^2(a_1b-1)\geq \frac{1}{4}a_1^2b^3+\frac{3}{4}b,$$
a contradiction with $c<0.25b^3$. 

Thus, $a_1d_1=a_2d_2$ must hold, so (\ref{relacija}) implies $(b+c-a_1-d_1)^2=(b+c-a_2-d_2)^2$ and since $a_i<b$ and $d_i<c$, we have $a_1+d_1=a_2+d_2$. Together with $a_1d_1=a_2d_2$, this yields $d_1=a_2$ and $d_2=a_1$, meaning that $\{a_1,a_2,b,c\}$ is a $D(4)$-quadruple.

\bigskip
\textbf{Acknowledgement:} The author was supported by the Croatian Science Foundation under the project no.\@ IP-2018-01-1313. Also, the author would like to thank Alan Filipin for his valuable comments on the earlier version of the manuscript.\\


\begin{thebibliography}{99}

\bibitem{aft} K. N. Adedji, A. Filipin, A. Togbé, {\it The problem of the extension of D(4)-triple \{1, b, c\}}, {  Rad Hrvat. Akad. Znan. Umjet. Mat. Znan.}, to appear.

\bibitem{bf} Lj.~Ba\'{c}i\'{c}, A.~Filipin, {\it On the extensibility of $D(4)$-pair $\{k-2,k+2\}$}, {  J.~Comb.~Number Theory} \textbf{5} (2013), 181--197.

\bibitem{mbt} M. Bliznac Trebješanin, {\it Extension of a Diophantine triple with the property $D(4)$}, {  Acta Math.~Hungar.} \textbf{163} (2021), 213–-246.

\bibitem{nas2} M.~Bliznac Trebje\v{s}anin, A.~Filipin, {\it Nonexistence of $D(4)$-quintuples}, {  J.~Number Theory}, \textbf{194} (2019), 170--217.

\bibitem{cdf} M. Cipu, A. Dujella and Y. Fujita, {\it Diophantine triples with largest two elements in common}, {  Period. Math. Hungar.} \textbf{82} (2021), 56--68.

\bibitem{duje_kon} A.~Dujella, {\it There are only finitely many Diophantine quintuples}, {  J.~Reine Angew.~ Math.~} \textbf{566} (2004), 183--214.

\bibitem{dujpet} A.~Dujella, A.~Peth\H{o}, {\it A generalization of a theorem of Baker and Davenport}, {  Quart.~J.~Math.~Oxford Ser.~(2)}, \textbf{49} (1998), 291--306.

\bibitem{dujram} A.~Dujella, A.~M.~S.~Ramasamy, {\it Fibonacci numbers and sets with the property $D(4)$}, {  Bull.~Belg.~Math.~Soc.~ Simon Stevin}, \textbf{12(3)} (2005), 401--412.

\bibitem{sestorka} A.~Filipin, {\it There does not exist a $D(4)$-sextuple}, {  J.~Number Theory} \textbf{128} (2008), 1555--1565.

\bibitem{fil_xy4} A.~Filipin,{ \it On the size of sets in which $xy + 4$ is always a square}, {  Rocky Mountain J.~ Math.~} \textbf{39} (2009), no.~ 4, 1195--1224.

\bibitem{fil_par} A.~Filipin, {\it The extension of some $D(4)$-pairs}, {  Notes Number Theory Discrete Math.~}\textbf{23} (2017), 126--135.

\bibitem{fht} A.~Filipin, Bo He, A.~Togb\'{e}, {\it  On a family of two-parametric D(4)-triples}, {  Glas.~Mat.~ Ser.~ III} \textbf{47} (2012), 31--51.

\bibitem{petorke} B.~He, A.~Togb\'{e}, V.~Ziegler, {\it There is no Diophantine quintuple}, {  Trans.~Amer.~Math.~Soc.~}\textbf{371} (2019), 6665--6709.

\bibitem{lm} M.~Laurent, M.~Mignotte, Yu.~Nesterenko, {\it Formes linéaires en deux logarithmes et déterminants d'interpolation}, {  J. Number Theory} \textbf{55} (1995), no. 2, 285--321. 

\bibitem{rickert} J.~H.~Rickert, {\it Simultaneous rational approximations and related Diophantine equations}, {  Proc.~Cambridge Philos.~Soc.~} \textbf{113} (1993) 461--472.


\end{thebibliography}
\end{document}